\documentclass[a4paper,10pt]{article}
\usepackage[cp1251]{inputenc}
\usepackage[english]{babel}

\usepackage{amssymb}\usepackage{amsmath}
\usepackage{epsfig}
\usepackage{graphics}
\def\bl{\rule[-1mm]{2.4mm}{2.4mm}}

\def\be{\begin{equation}}
\def\ee{\end{equation}}

\newtheorem{thrm}{\bf Theorem}
\newtheorem{lmm}{\bf Lemma}

\begin{document}

\title {Image of Abel-Jacobi map\\ for hyperelliptic genus 3 and 4 curves}
\author{\copyright 2013 ~~~~A.B.Bogatyrev
\thanks{Supported by RFBR grants 13-01-00115 and RAS Program
"Modern problems of theoretical mathematics"}}
\date{}
\maketitle

\hfill\parbox{8cm}{\it To the memory of Andrei Gonchar and Herbert Stahl}
\vspace{1mm}

The necessity for computation of abelian integrals often arises in problems of classical mechanics (see \cite{Go,D} and references therein), conformal mappings of polygons \cite{B,G},  
solitonic dynamics, general relativity \cite{EHKKL}, solutions of rational optimization problems \cite{B2} and so on. Function theory on Riemann surfaces allows one to evaluate (with computer accuracy) such integrals without any quadrature rules. As an example, let us consider Riemann's formula for the abelian integral with two simple poles at the points $R,Q$ of the curve $\cal X$:

\be
\eta_{RQ}(P):=\int^P d\eta_{RQ}=
\log\frac{\theta[\epsilon,\epsilon'](u(P)-u(R),\Pi)}
{\theta[\epsilon,\epsilon'](u(P)-u(Q),\Pi)}+const,
\label{intrep}
\ee
where $\theta[\cdot](\cdot,\cdot)$ is Riemann theta function with some odd integer characteristics $[\epsilon,\epsilon']$. 
It might seem that the usage of this formula 
still requires the computation of holomorphic abelian integrals involved in the Abel-Jacobi (AJ) map:

\be
\label{AJmap}
u(P):=\int_{P_0}^P du\in\mathbb{C}^g,
\qquad du:=(du_1,du_2,\dots,du_g)^t,
\ee 
where $g$ is the genus of the curve $\cal X$ and $du_s$ are suitably normalized abelian differentials of the first kind.
In this note we show how to avoid the evaluation of AJ map: the image of low genus $g<5$ hyperelliptic curve in its Jacobian
is given as the solution of a (slightly overdetermined) set of equations including theta functions. 
Moving along the curve embedded in its Jacobian we can compute abelian integral  by explicit formula (e.g. (\ref{intrep})) 
and simultaneously compute some projection of the curve to the complex projective line (see e.g. \cite{FK,B}). In this way we get a parametric representation of abelian integrals \cite{B,G} which may be used either for the evaluation of integral, or for its inversion. See also \cite{EHKKL} for the alternative approach using 
hyperelliptic $\sigma$ and $\wp$ functions.

\section{Introduction}
\begin{figure}[h!]
\begin{picture}(150,55)
\put(0,0){\includegraphics[scale =.75]{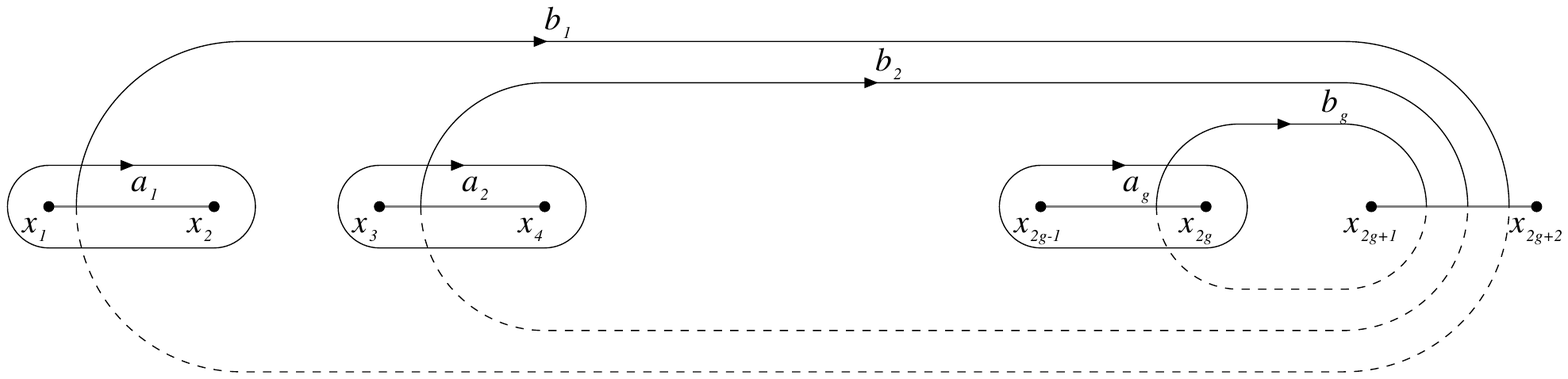}}
\end{picture}
\caption{\small Symplectic basis in the homologies of curve $X$}
\label{Basis}
\end{figure}

Let us fix the notations. Consider genus $g$ hyperelliptic curve $\cal X$:  
$$
w^2=\prod_{j=1}^{2g+2}(x-x_j),
$$
with distinct branch points $x_1,x_2,\dots,x_{2g+2}$. The curve admits involution $J(x,w):=(x,-w)$
with fixed points $P_s=(x_s,0)$, $s=1,\dots,2g+2$.
We introduce a symplectic basis in the homologies of $\cal X$ as in the Fig. \ref{Basis}.
Dual basis of holomorphic differentials satisfies normalization
$$
\int_{a_j}du_s:=\delta_{js},
$$
and generates the period matrix
$$
\int_{b_j}du_s=:\Pi_{js}.
$$

It is convenient to represent points  $u\in\mathbb{C}^g$ as theta characteristics, i.e. couple of real g-vector columns $\epsilon, \epsilon'$:
\be
u=(\epsilon'+\Pi\epsilon)/2.
\ee
The points of Jacobian  
\be
Jac({\cal X}):=\mathbb{C}^g/L(\Pi), 
\qquad L(\Pi):=\mathbb{Z}^g+\Pi\mathbb{Z}^g,
\label{Jac}
\ee
in this notation correspond to two vectors with real entries modulo 2. Second order points of Jacobian are
$2\times g$ matrices with binary entries. In particular, the 
the images of the Weierstrass points of $\cal X$ under the AJ map (\ref{AJmap}) with the initial point $P_0:=P_{2g+2}$ 
are as follows:

\begin{tabular}{c|c|c}
$P_s$ & $u(P_s)~ mod~ L(\Pi)$ & $[\epsilon, \epsilon']^t$\\
\hline
$P_1$ & $\Pi_1/2$ & $\tiny
\left[\begin{array}{c} 100\dots\\000\dots\end{array}\right]$\\
$P_2$&$(\Pi_1+E_1)/2$&$\tiny
\left[\begin{array}{c} 100\dots\\100\dots\end{array}\right]$\\
$P_3$&$(\Pi_2+E_1)/2$&$\tiny
\left[\begin{array}{c} 010\dots\\100\dots\end{array}\right]$\\
$P_4$&$(\Pi_2+E_1+E_2)/2$&$\tiny
\left[\begin{array}{c}010\dots\\110\dots\end{array}\right]$\\
$P_5$&$(\Pi_3+E_1+E_2)/2$&$\tiny
\left[\begin{array}{c}0010\dots\\1100\dots\end{array}\right]$\\
\vdots&\vdots&\vdots\\
$P_{2g-1}$&$(\Pi_g+E_1+E_2+\dots+E_{g-1})/2$&$\tiny
\left[\begin{array}{c} 0\dots01\\1\dots10\end{array}\right]$\\
$P_{2g}$&$(\Pi_g+E_1+E_2+\dots+E_g)/2$&$\tiny
\left[\begin{array}{c} 0\dots01\\11\dots1\end{array}\right]$\\
$P_{2g+1}$&$(E_1+E_2+\dots+E_g)/2$&$\tiny
\left[\begin{array}{c}00\dots0\\11\dots1\end{array}\right]$\\
$P_{2g+2}$&$0$&$\tiny
\left[\begin{array}{c}00\dots0\\00\dots0\end{array}\right]$
\label{AJPj}
\end{tabular}   

Here $E_s$ and $\Pi_s$ are the columns of the identity matrix and the period matrix respectively.

The following series has very high convergence rate and well controlled accuracy \cite{DHB}.
\be
\label{thetadef}
\theta(u, \Pi):=\sum\limits_{m\in\mathbb{Z}^g}
\exp(2\pi m^tu+\pi i m^t\Pi m),
\qquad u\in\mathbb{C}^g;\Pi=\Pi^t\in\mathbb{C}^{g\times g};\Im\Pi>0.
\ee
Matrix argument $\Pi$ of theta function may be omitted if it does not lead to a confusion.
This function (\ref{thetadef}) has the following easily checked quasi-periodicity properties with respect to the lattice $L(\Pi)$:
\be
\label{quasiperiod}
\theta(u+m'+\Pi m; \Pi)=\exp(-i\pi m^t\Pi m-2i\pi m^tu)\theta(u;\Pi),
\qquad m,m'\in\mathbb{Z}^g.
\ee

Theta function may be considered as a multivalued function in the Jacobian
or as a section of a certain line bundle. The zero set of theta function -- the theta divisor -- is well defined 
in the Jacobian since the  factors in right hand side of (\ref{quasiperiod}) do not vanish.
The theta divisor is described by so called Riemann vanishing theorems  \cite{RF,FK}.
One of the important ingredients in those theorem is a vector of Riemann's constants $\cal K$
which depends on the choice of homology basis and the initial point in AJ map. 
In the above setting the vector of Riemann's constants may be found by a straightforward computation \cite{Mu} or by 
some combinatorial argument \cite{FK} and corresponds to a characteristic 
$$
{\cal K}(P_0) \sim \left[\begin{array}{c}\dots11111\\ \dots10101\end{array}\right]^t.
$$

Theta divisor is described by the following
\begin{thrm}[Riemann] $\quad\theta(e)=0$ iff $e=u(D_{g-1})+{\cal K}\quad mod~L(\Pi)$
for some degree $g-1$ positive divisor $D_{g-1}$ on the curve.
\label{ThDiv}
\end{thrm}
In particular this means that genus 2 curve in its Jacobian is just the solution of one equation
$\theta(u+{\cal K})=0$. Now we consider higher genera.

\section{Genus three hypereliptic curves}
\begin{thrm} 
Let $P,Q$ be two distinct points on the genus 3 hyperelliptic curve $\cal X$.
The solution of the set of two equations 
\be
\label{theta2}
\begin{array}{c}
\theta(u-u(P)+{\cal K})=0,\\
\theta(u-u(Q)+{\cal K})=0,
\end{array}
\ee
in the Jacobian $Jac(X)$ is the union of $u({\cal X})$ -- the image of the curve under AJ map -- 
and its shift $\quad u({\cal X})+u(P+Q)$. 
\end{thrm}
{\bf Proof.} Let $u$ be a point in the intersection of two shifted theta divisors (\ref{theta2}).
The representation of the theta divisor from Theorem \ref{ThDiv}, suggests that there are two degree 2 positive divisors $D_2$ and $D'_2$ on the curve $\cal X$ such that 
$$
u=u(D_2+P)=u(D'_2+Q) \quad mod~L(\Pi).
$$
We consider two cases:\\
1) When two linearly equivalent divisors $D_2+P\sim$ $D'_2+Q$ are non-special,
they coincide and therefore
$
u\in u({\cal X})+u(P+Q)
$ 
in the Jacobian of the curve.\\
2) When $i(D_2+P)=$ $i(D'_2+Q)>0$,
each divisor contains $J$ -equivalent points which together give zero input to the AJ map.
Therefore $u\in u({\cal X})$. 

Conversely, consider any point $u$ in  the union of $u({\cal X})$ and its shift by $u(P+Q)$.
In other words, $u=u(S)$ or $u=u(S+P+Q)$ for some point $S$ of our curve. The first equation in the system
(\ref{theta2}) is true because the argument of theta function satisfies the condition of Theorem \ref{ThDiv}
with $D_2:=S+JP$ or $D_2:=S+Q$ respectively. Same argument (with proper choice of $D_2$ ) fits for the other equation. 
~~~\bl

We see that two equations are not enough to localize the image of AJ map since the parasitic component 
$u({\cal X})+u(P+Q)$ arises. Adding yet another equation of this type, we achieve the goal.

\begin{thrm} 
Let $P,Q,R$ be two distinct points on the genus 3 hyperelliptic curve $\cal X$.
The solution of the set of three equations 
\be
\label{theta3}
\begin{array}{c}
\theta(u-u(P)+{\cal K})=0,\\
\theta(u-u(Q)+{\cal K})=0,\\
\theta(u-u(R)+{\cal K})=0
\end{array}
\ee
in the Jacobian $Jac({\cal X})$ is the union of $u({\cal X})$ -- the image of the curve under AJ map -- 
and just one point $u(P+Q+R)$. 
\end{thrm}
{\bf Proof.} Taking into account the previous theorem, we have to find the intersection of two different shifts of the same curve 
$u({\cal X})$: by $u(P+Q)$ and by $u(P+R)$. The point in this intersection has two representations
$$
u=u(S+P+Q)=u(S'+P+R)
$$
for some points $S,S'$ of the curve. If $i(S+P+Q)=i(S'+P+R)=0$
then two mentioned divisors coincide: $S=R$, $S'=Q$ and therefore $u=u(P+Q+R)$.
Assume that speciality index is positive. Then each of the divisors contains
two J-equivalent points whose AJ-images are opposite. Hence, $u\in u({\cal X})$. ~~~\bl

{\bf Remark.} 
There is a temptation to put $P=JQ$ in the system (\ref{theta2}) and to get the curve $u(\cal X)$ as a solution
of just two equations. However, the only points of the curve with known value of AJ map are its 
Weierstrass points $P_s$ -- see the table above. With this choice of auxiliary points, the systems (\ref{theta3}) and the forthcoming system (\ref{theta4}) may be rewritten in terms of theta functions with integer characteristics.

\section{Genus four hypereliptic curves}
We need an auxiliary statement:
\begin{lmm}
\label{Jpts}
Let $D_s$ be a degree $s$ positive divisor on genus 4 hyperelliptic curve and $u(D_3)=u(D_5)$. Then $D_5$ contains $J$-equivalent points. 
\end{lmm}
{\bf Proof} Let us detach a point from the larger divisor: $D_5=D_4+P$, then
$D_4$ is linear equivalent to $D_3+JP$. Now if $i(D_4)>0$ then $D_4$ 
contains $J$-equivalent points. Otherwise $D_4$ is non-special and it contains  
$JP$, hence $D_5\ge P+JP$. ~~\bl

\begin{thrm} 
Let $P_1,P_2$; $Q_1,Q_2$ be four distinct points on the genus 4 hyperelliptic curve $\cal X$
and  $P_1\neq JP_2$; $Q_1\neq JQ_2$. The solution of the set of three equations 
\be
\label{theta3*}
\begin{array}{r}
\theta(u+{\cal K})=0,\\
\theta(u-u(P_1+P_2)+{\cal K})=0,\\
\theta(u-u(Q_1+Q_2)+{\cal K})=0,
\end{array}
\ee
in the Jacobian is the union of $u(\cal X)$ -- the image of the curve under AJ map -- 
and its shifts by four vectors $u(P_j+Q_s)$, $j,s=1,2$. 
\end{thrm}

{\bf Remark} The condition of this theorem implies that $u(P_1+P_2)\neq u(Q_1+Q_2) ~~mod~L(\Pi)$.
Indeed, if it were not the case, then by Abel's theorem $P_1$ and $P_2$ were zeroes of degree 2
function from $\mathbb{C}({\cal X})$. This function is essentially unique and therefore $P_1=JP_2$.

{\bf Proof}. Let $u$ be a point in the intersection of three shifted theta divisors (\ref{theta3*}).
Due to the representation of the theta divisor, there are three positive divisors $D_3$, $D'_3$, $D''_3$,
each of degree 3 such that 
$$
u=u(D_3)=u(D'_3+P_1+P_2)=u(D''_3+Q_1+Q_2)
\quad mod~L(\Pi).
$$
The divisors $D'_5:=D'_3+P_1+P_2$ and $D''_5:=D''_3+Q_1+Q_2$ in two latter equations contain $J$-equivalent points in accordance with the Lemma \ref{Jpts}.
We consider two cases:\\
1) When $i(D_3)>1$, then the divisor contains two points which give opposite input to the AJ mapping.
Hence, $u\in u({\cal X})$.
2) When $i(D_3)=1$ (non-special divisor), both divisors $D'_5$ and $D''_5$ with $J$-equivalent points thrown away,
coincide with $D_3$. Now $D_3$ contains at least one point of $P_1, P_2$ and one point of $Q_1,Q_2$.
Hence, $u\in u({\cal X})+u(P_j+Q_s)$ for some $j$ and $s$. 

Conversely, for the points $u$ in the shifted AJ-images of the curve $\cal X$, the arguments of the theta functions in 
(\ref{theta3*}) satisfy the condition of the Theorem \ref{ThDiv}. Say, for the point $u=u(S+P_1+Q_2)$,  $S\in{\cal X}$,
the divisor $D_3=S+P_1+Q_2$ for the first equation; $D_3=S+JP_2+Q_2$ for the second equation; $D_3=S+P_1+JQ_1$ for the third equation in (\ref{theta3*}). For $u=u(S)$, the divisor $D_3=S+P_1+JP_1$ for the first equation; $D_3=S+JP_1+JP_2$ for the second equation; $D_3=S+JQ_1+JQ_2$ for the third equation in the system.
~~~\bl 

Adding yet another equation of this type to (\ref{theta3*}) allows us to localize the image of AJ map for genus 4 hyperelliptic curves.  Eight parasitic points however arise in the solution. 

\begin{thrm} 
Let $P_s,Q_s,R_s$, $s=1,2$ be six distinct points on the genus 4 hyperelliptic curve $\cal X$
and the points of each of three pairs $P_s,Q_s,R_s$ are not $J$-equivalent. The solution of the set of four equations 
\be
\label{theta4}
\begin{array}{r}
\theta(u+{\cal K})=0,\\
\theta(u-u(P_1+P_2)+{\cal K})=0,\\
\theta(u-u(Q_1+Q_2)+{\cal K})=0,\\
\theta(u-u(R_1+R_2)+{\cal K})=0,
\end{array}
\ee
in the Jacobian $Jac({\cal X})$ is the union of $u({\cal X})$ -- the image of the curve under AJ map -- 
and eight points $u(P_j+Q_k+R_s)$, ~$j,k,s=1,2$. 
\end{thrm}
{\bf Proof.} Again, we have to find the intersection of the solutions obtained in the previous theorem.
Let us shift the same curve 
$u({\cal X})$ by $u(P_s+Q_j)$ and by $u(P_m+R_k)$. The point $u$ in this intersection has two representations
$$
u=u(S+P_s+Q_j)=u(S'+P_m+R_k)
$$
for some points $S,S'$ of the curve. If $i(S+P_s+Q_j)=i(S'+P_m+R_k)=1$
then two mentioned non-special divisors coincide: this may only happen in case $m=s$, $S=R_k$ and $S'=Q_j$ and therefore $u=u(P_s+Q_j+R_k)$.
Assuming that speciality index is greater than 1, we come to a conclusion that $u\in u({\cal X})$. ~~~\bl

\vspace{5mm}
\parbox{7cm}
{\it
Institute for Numerical Mathematics,\\
Russian Academy of Sciences;\\
Moscow Inst. of Physics and Technology;
Moscow State University\\[1mm]
{\tt gourmet@inm.ras.ru, ab.bogatyrev@gmail.com}}

\end{document}